\documentstyle[12pt]{article}
\begin{document}

\begin{flushright}
MRC-PH-TH-11-96

q-alg/9609028
\end{flushright}

\begin{center}

{\bf \Large A q--oscillator Green Function} 

\vspace{1cm}

H. Ahmedov$\ast$  and I.H.Duru$\ast\dag$ 
\end{center}
$\ast$ -TUBITAK -Marmara Research Centre, Research Institute for Basic
Siences, Department of Physics, P.O. Box 21, 41470 Gebze, Turkey

$\dag$ -Trakya University, Mathematics Department, P.O. Box 126, Edirne, Turkey.

\vspace{1cm}

\begin {center}
Abstract
\end{center}

By using the generating function formula for the product of two q-Hermite
polynomials q-deformation of the Feynman Green function for the harmonic 
oscillator is obtained.
\

\

\

PACS numbers: 03.65.Fd and 02.20.- a
\begin {center}
September 1996
\end{center}

\vfill
\eject

\section{Introduction}
q-oscillators are the most extensively studied deformed dynamical
systems. They have been presented in severent different types \cite{kn:a}.
The e-functions of these q-oscillators are expressible in terms of either
the discrete or the continuous q-Hermite polynomials \cite{kn:v}.

Although the literature on q-oscillators is very rich, the corresponding 
q-Green functions have not been investigated. This looks suprising when one 
considers the fact that the exact closed form of the Green function of the undeformed
oscillator has been known for many decades \cite{kn:fey}.

It is the purpose of this note to obtain the Green function for one of the 
q-oscillators. The q-oscillator we deal with is the one which is solved
in terms of the continuous $q^{-1}$-Hermite polynomials \cite{kn:ata}.

In Section II we briefly review the q-oscillator realization of Ref.4.
Section III is devoted to the derivation of the q-oscillator Green function
which is the deformation of the well known Feynman formula. The method of the 
calculation of the non-trivial $q\rightarrow 1$ limit , which is essential
for arriving at the usual Feynman Green function is outlined in the Appendix.

\section{A q-oscillator realization}
Recently Atakishiev, Frank and Wolf introduced a simple difference 
realization of the Heysenberg q-algebra \cite{kn:ata}.They also studied 
the corresponding q-oscillator Hamiltonian and its e-functions in terms
of the $q^{-1}$-Hermite polynomials.The q-annihilation and creation 
operators acting on the smooth functions $f(\xi)$ with $\xi\in(-\infty,\infty)$
are given by

\begin{equation}
b_q=
\frac{1}{2^{1/2}q^{1/4}\beta}
v(\xi)(q^{-\beta\xi} \exp(\frac{1}{2\beta}\partial_{\xi})  -
q^{\beta\xi} \exp(-\frac{1}{2\beta}\partial_{\xi})) v(\xi)
\end{equation}

\begin{equation}
b^{\dag}_q=
\frac{1}{2^{1/2}q^{1/4}\beta}
v(\xi)(q^{-\beta\xi} \exp(-\frac{1}{2\beta}\partial_{\xi}) - 
q^{\beta\xi} \exp(\frac{1}{2\beta}\partial_{\xi})) v(\xi)
\end{equation}
where

\begin{equation}
k = -\log q, \ \ \  \beta = \frac{1}{(2(1-q))^{1/2}}, \ \ \ v(\xi)=\frac{1}{(\cosh(k\beta\xi))^{1/2} }   
\end{equation}
and $\xi$ is the dimensionless variable (with $\hbar$ =1)

\begin{equation}
\xi=\sqrt{\omega m} x.
\end{equation}          
The algebra satisfied by the operators (1), (2) is

\begin{equation}
b_q b^{\dag}_q - qb^{\dag}_q b_q  =1.
\end{equation}
In the limit $q\rightarrow 1^{-} (k\rightarrow 0^{+})$ 
$b_q$, $b^{\dag}_q$ takes the usual forms:

\begin{equation}
b=\frac{1}{2^{1/2}}(\xi+\frac{d}{d\xi})  \ \ \ and \ \ \
b^{\dag}=\frac{1}{2^{1/2}}(\xi-\frac{d}{d\xi}). 
\end{equation}
The operator

\begin{equation}
H_q = b^{\dag}_q b_q
\end{equation}
which is self-adjoint under the inner product 

\begin{equation}
(f,g)=\int^{\infty}_{-\infty} d\xi 
\overline{f(\xi)}g(\xi)
\end{equation}
satisfies the e-value equation

\begin{equation}
H_q \Psi_n^q (\xi) = [n] \Psi_n^q (\xi). 
\end{equation}
Here [n] is defined as usual as

\begin{equation}
[n] = \frac{1-q^n}{1-q},  \ \ \ \ n =0,1,2 ...
\end{equation}
and the e-functions are given by

\begin{equation}
\Psi_n^q (\xi) = (\frac{k}{\pi(1-q)})^{1/4}
\frac{q^{(n+1/2)^2/4}} {((q;q)_n)^{1/2}}
(\cosh(k\beta\xi))^{1/2} \exp(-k\beta^2\xi^2) h_n(\sinh(k\beta\xi)\mid q)
\end{equation}
$h_n$ is the continuous $q^{-1}$- Hermite polynomial
 and $(q;q)_n$ is the q- factorial :

\begin{equation}
(q;q)_n =\prod^n_{j=1} (1-q^j)
\end{equation}

In $q\rightarrow 1^{-}$ limit $h_n$ takes the form of the usual Hermite 
polynomial (with  $\sinh(k\beta\xi)\rightarrow (\frac{1-q}{2})^{1/2}\xi$):

\begin{equation}                                                        
\lim_{q\rightarrow 1^{-}} (\frac{2}{1-q})^{n/2} h_n(\sinh(k\beta\xi)\mid q) =
H_n(\xi).
\end{equation}

\section{A ``physical" q-oscillator and its Green function}

Making use of the operator (7) we can write the following ``physical"
q- oscillator  Schr\"{o}dinger equation including the ground state energy :

\begin{equation}
(q^{1/2}b^{\dag}_q b_q + [\frac{1}{2}]) \Phi_n^q (\xi ,t) = 
\zeta D^q_{\zeta} \Phi_n^q (\xi ,t).
\end{equation}
Here $\zeta$ is the exponantial time parameter given by 

\begin{equation}
\zeta = \exp(-i\omega t)
\end{equation}
and $\Phi_n^q (\xi ,t)$ is the time dependent wave function :
 
\begin{equation}
\Phi_n^q (\xi ,t) =\exp(-i\omega (n+1/2)) \Psi_n^q (\xi ) = 
\zeta^{n+1/2} \Psi_n^q (\xi ) .
\end{equation}

The action of the q- derivative on the time dependent factor of the above 
wave function

\begin{equation}
\zeta D^q_{\zeta} \zeta^{n+1/2}=[n+1/2]\zeta^{n+1/2} =  
(q^{1/2}[n]+[\frac{1}{2}]) \zeta^{n+1/2} 
\end{equation}
exhibits the correct energy spectrum of the "physical`` q- oscillator.

\

\

Time dependent wave function enable us to write the q- Green function
for the oscillator as

\begin{equation}
K_q(\xi,\xi';z)=
\sum_{n=0}^{\infty} z^{n+1/2}\Psi_n^q (\xi ) \Psi_n^q (\xi')
\end{equation}
where

\begin{equation}
z=\overline{\zeta}\zeta'=\exp(-i\omega (t'-t)).
\end{equation}

To execute the summation over n in (18) we recall the following generating 
function formula for the product of two continuous q- Hermite polynomials 
\cite{kn:vil} :

\begin{eqnarray}
\frac{(z^2;q)_{\infty}}
{(z\exp(i(\theta +\phi));q)_{\infty} 
(z\exp(i(\theta -\phi));q)_{\infty} (z\exp(-i(\theta +\phi));q)_{\infty}}\times  \nonumber \\
\times\frac{1}{(z\exp(-i(\theta -\phi));q)_{\infty}}
=\sum_{n=0}^{\infty}z^n 
\frac{H_n(\cos\theta\mid q) H_n(\cos\phi\mid q)}{(q;q)_n}
\end{eqnarray}
Here $(\alpha;q)_{\infty}$ is defined as

\begin{equation}
(\alpha;q)_{\infty}=\prod_{j=0}^{\infty} (1-\alpha q^j)
\end{equation}

In $q\rightarrow 1^{-}$ limit (see Appendix) the formula (20) is reduced to the well known 
summation formula for the product of two undeformed Hermite polynomials 
\cite{kn:mor}:

\begin{eqnarray}
\frac{1}{(1-z^2)^{1/2}}
\exp[-\frac{1}{1-z^2}( z^2(\xi^2 + \xi'^2)-2z\xi\xi')] = \nonumber \\
\sum_{n=0}^{\infty}
\frac{z^n}{2^n n}H_n(\xi)H_n(\xi')
\end{eqnarray}
with

\begin{equation}
\xi=(\frac{1-q}{2})^{1/2}\cos\theta,  \ \ \ \ \ \ \ 
\xi'= (\frac{1-q}{2})^{1/2}\cos\phi.
\end{equation}
Note that by the help of (22) one can derive the well known Feynman formula 
for the undeformed oscillator \cite{kn:fey} (with $T=t'-t$)

\begin{equation}
K(\xi ,\xi' ;t'-t)=(\frac{m\omega}{2\pi i \sin(\omega T)})^{1/2}
\exp[\frac{im\omega}
{2\sin(\omega T)}
((\xi^2 +\xi'^2)\cos(\omega T) - 2\xi\xi')]
\end{equation}
from the Green function written in the wave function 
decomposition form \cite{kn:dur}.

\

To derive the q- oscillator Green function we first insert $q^{-1}$ in place of q
by recalling the relation \cite{kn:vil}

\begin{equation}
\frac{1}{(q^{-1};q^{-1})_n } =
q^{1/8} q^{(n+1/2)^2/2} \frac{(-1)^n}{(q;q)_n}.
\end{equation}
After making the required analytic continuations in $\theta$ and $\phi$
we arrive at

\begin{eqnarray}
E^{-1}_q(\frac{qz^2}{1-q})
E_q(\frac{-qz}{1-q}\exp-(\theta +\phi))
E_q(\frac{-qz}{1-q}\exp(\theta +\phi)) \nonumber \\ 
E_q(\frac{qz}{1-q}\exp-(\theta -\phi))
E_q(\frac{qz}{1-q}\exp(\theta -\phi))  \nonumber \\
=\sum_{n=0}^{\infty}
\frac{h_n(\sinh\theta\mid q) h_n(\sinh\phi\mid q)}
{(q;q)_n}q^{(n+1/2)^2/2}  z^n
\end{eqnarray}
which is the formula suitable to our $q^{-1}$- Hermite polynomials.
The q- exponantials employed in the above equation are qiven in terms  
of the $n\rightarrow\infty$ limit of the q- factorials as \cite{kn:vil} 

\begin{equation}
E_{1/q}(-x)=E_q^{-1}(x)=((1-q)x;q)_{\infty}.
\end{equation}

When we introduce the formula of (26) into (18) we obtain the final form
of the q-oscillator Green function:

\begin{eqnarray}
K_q(\xi ,\xi' ;z)=  
q^{-1/8}(\frac{2k}{\pi } )^{1/2} \beta 
(\cosh(k\beta\xi ) \cosh(k\beta\xi' ))^{1/2}
\exp(-k\beta^2(\xi^2 +\xi'^2) z^{1/2} \nonumber \\
E_q^{-1}(\frac{qz^2}{1-q})
E_q(\frac{-qz}{1-q}\exp(-k\beta (\xi +\xi')) 
E_q(\frac{-qz}{1-q}\exp(k\beta (\xi +\xi'))  \nonumber \\
E_q(\frac{qz}{1-q}\exp(-k\beta (\xi -\xi')) 
E_q(\frac{qz}{1-q}\exp(k\beta (\xi -\xi')) 
\end{eqnarray}

By the process sketched in the Appendix the above equation is 
reduced to the Feynman formula of (24) in $q\rightarrow 1^{-}$ limit.

In T$\rightarrow 0$ (z$\rightarrow 1$) limit we distinguish two cases :

(i) For $\xi\neq\xi'$  by the virtue of the first exponantial
\begin{equation}
E_q^{-1}(\frac{qz^2}{1-q})=(qz^2;q)_{\infty}=
\prod_{n=1}^{\infty}(1-q^nz^2).
\end{equation}
we have
\begin{equation}
\lim_{T\rightarrow 0}K_q(\xi ,\xi' ;z)=0.
\end{equation}

(ii) For $\xi = \xi'$ on the other hand by the virtue of the
$1^{\underline{st}}$, $2^{\underline{th}}$ and $3^{\underline{th}}$ expanantiols we
have
\begin{eqnarray}
\lim_{z\rightarrow 1}
E_q^{-1}(\frac{qz^2}{1-q})
E_q(\frac{-qz}{1-q}) 
E_q(\frac{-qz}{1-q}) \nonumber \\
=\lim_{z\rightarrow 1}
\prod_{n=0}^{\infty}\frac{(1-q^nz^2)}{(1-q^nz)(1-q^nz)}.
\end{eqnarray}
The n=1 factor in the above equation contributes a singularity of
\begin{equation}
\lim_{z\rightarrow 1}\frac{1}{1-z}
\end{equation}
type.

It is easy to conclude then that the Green function (28) behaves as the 
$\delta$-function $\delta(\exp(k\beta\xi)-\exp(k\beta\xi`))$ in
 T$\rightarrow 0$ (z$\rightarrow 1$) limit. 

\

Acknowledgement. We thank \"{O}.F.Dayi for reading the manuscript. 

\

\begin{center}
APPENDIX
\end{center}

\renewcommand{\theequation}{A.\arabic{equation}}
\setcounter{equation}{0}

Using the definitions in (21) and (23) we can rewrite (20) as

\begin{eqnarray}
\prod_{j=0}^{\infty} (1-z^2q^j)(1-z^2q^{2j})^{-2}(1-(1-q)A_j)^{-1}= \nonumber \\
\sum_{n=0}^{\infty} \frac{z^n}{(q;q)_n}
H_n((\frac{1-q}{2})^{1/2}\xi\mid q) H_n((\frac{1-q}{2})^{1/2}\xi'\mid q)
\end{eqnarray}
with

\begin{equation}
A_j =2\frac
{ zq^j\xi\xi'(1+z^2q^{2j})-(\xi^2 +\xi'^2)z^2q^{2j} }
{(1-z^2q^{2j})^2}
\end{equation}
$q\rightarrow 1^{-}$ limit of the right hand side (r.h.s) of (A.1)
( with $(q;q)_n=(1-q)^n [n] $) is

\begin{equation}
\lim_{q\rightarrow 1^{-}}(r.h.s) =
\sum_{n=0}^{\infty} \frac{z^n}{2^n n} H_n(\xi) H_n(\xi')
\end{equation}

Let us take the logarithm of the left hand side (l.h.s.) of (A.1) :

\begin{equation}
\log (l.h.s.) =\sum_{j=0}^{\infty} [ \log (1-z^2q^j)-2\log (1-z^2q^{2j}) ] \nonumber \\
- \sum_{j=0}^{\infty} \log (1-(1-q)A_j)
\end{equation}

Expanding the logarithm function into the series (for $\mid z\mid\leq 1$  )
the first two terms in the above equation can be written as

\begin{eqnarray}
\sum_{j=0}^{\infty} [ \log (1-z^2q^j)-2\log (1-z^2q^{2j}) ] =
\sum_{k=1}^{\infty} \frac{z^{2k}}{k(1+q^k)} \nonumber \\
=-\log (1-z^2)_{q^2} ^{1/2}
\end{eqnarray}
   
with \cite{kn:vil}
 
\begin{equation}
(1-z)_q^a= \psi_{1,0}(q^{-a},q;z)
\end{equation}

Then (l.h.s.) can be rewritten as

\begin{equation}
(l.h.s.)=\frac{1}{(1-z^2)_{q^2}^{1/2}}\exp(- F(\xi ,\xi',z\mid q))
\end{equation}

where

\begin{equation}
F(\xi ,\xi',z\mid q)=\sum_{j=0}^{\infty}\log (1-(1-q)A_j) 
\end{equation}
>From (A.2) we see that the functional sequence  $\mid A_j \mid $  for any 
value of $z$, $\xi$ and $\xi'$ (except the case z=1)  decreases :
\begin{equation}
\mid A_0 \mid \ > \mid A_1 \mid \ ... > \mid A_n \mid \ ...                
\end{equation}                                      
Let us consider the zeroth term of the sequence (A.9) $\mid A_0 \mid$ 
and fix the value of z (with z$\neq$1). It is clear that there exist
$n\in N$ such that

\begin{equation}
(1-q_i)\mid A_0 \mid <1 \ \ \  \ \  i>n.
\end{equation}
where  ${q_i}$ is the sequence: $\lim_{i\rightarrow \infty}q_i=1^{-}$.
By the virtue of (A.9) the functions   $(1-q_i)\mid A_j \mid $ satisfy  the 
condition (A.10) too.Thus the logarithm function (A.9) 
can be expanded in the Taylor series in $q\rightarrow 1^{-}$ 
limit as

\begin{equation}
\lim_{q\rightarrow 1^{-}}F(\xi ,\xi',z\mid q)=
\lim_{q\rightarrow 1^{-}}
\sum_{k=1}^{\infty}\frac{(1-q)^k}{k}
\sum_{j=0}^{\infty}(A_j)^k 
\end{equation}
In the above expression only the k=1 term survives:
\begin{equation}
\lim_{q\rightarrow 1^{-}}F(\xi ,\xi',z\mid q)=
\lim_{q\rightarrow 1^{-}}
(1-q)
\sum_{j=0}^{\infty}A_j 
\end{equation}
After expanding the denominator of $A_j$ into the power 
serries we arrive at                    

\begin{equation}
\lim_{q\rightarrow 1^{-}} F(\xi ,\xi',z\mid q) =
\frac{z^2(\xi^2 +\xi'^2) -2z\xi\xi'}{1-z^2}
\end{equation}

>From (A.8) and (A.14) we get

\begin{equation}
\lim_{q\rightarrow 1^{-}}(l.h.s.)=\frac{1}{(1-z^2)^{1/2}}
\exp[-\frac{z^2(\xi^2 +\xi'^2) -2z\xi\xi'}{1-z^2}] 
\end{equation}

which together with (A.3) establishes the desired limit of (22).

\vfill
\eject

\end{document}